\documentclass[%
 aip,
 amsmath,amssymb,
 reprint,%
]{revtex4-1}

\usepackage{graphicx}
\usepackage{dcolumn}
\usepackage{bm}

\usepackage[utf8]{inputenc}
\usepackage[T1]{fontenc}
\usepackage{mathptmx}
\usepackage{etoolbox}
\usepackage{color}
\usepackage{amsthm, amssymb}

\theoremstyle{plain}

\newtheorem{thm}{{Theorem}}
\newtheorem*{thm-non}{{Theorem}} 
\newtheorem{lemma}{{Lemma}}

\newtheorem*{defi}{{Definition}}
\newtheorem*{remark}{\emph{Remark}}

\newtheorem{pr*}{Problem*}

\newtheorem*{ex-pr}{Extra Problem}

\makeatletter
\def\@email#1#2{%
 \endgroup
 \patchcmd{\titleblock@produce}
  {\frontmatter@RRAPformat}
  {\frontmatter@RRAPformat{\produce@RRAP{*#1\href{mailto:#2}{#2}}}\frontmatter@RRAPformat}
  {}{}
}%
\makeatother
\begin{document}

\preprint{AIP/123-QED}

\title[Symmetry Groupoids for Pattern-Selective Feedback Stabilization of the Chafee--Infante Equation]{Symmetry Groupoids for Pattern-Selective Feedback Stabilization of the Chafee--Infante Equation
}
\author{I. Schneider}
 \altaffiliation[Also at ]{Freie Universität Berlin, Institut für Mathematik, Arnimallee 7, 14195 Berlin, Germany}

\affiliation{ 
Universität Rostock, Institut für Mathematik, Ulmenstr. 69, 18057 Rostock, Germany
}%

\author{J.-Y. Dai}
 \email{jydai@nchu.edu.tw}
\affiliation{%
National Chung Hsing University, Department of Applied Mathematics, 145 Xingda Rd., Taichung City, Taiwan
}%

\date{\today}

\begin{abstract}
Reaction-diffusion equations are ubiquitous in various scientific domains and their patterns represent a fascinating area of investigation. However, many of these patterns are unstable and therefore challenging to observe. 
To overcome this limitation, we present new noninvasive feedback controls based on symmetry groupoids. As a concrete example, we employ these controls to selectively stabilize unstable
equilibria of the Chafee--Infante equation under Dirichlet boundary conditions on the interval. Unlike conventional reflection-based control schemes, our approach incorporates additional symmetries that enable us to design new convolution controls for stabilization. By demonstrating the efficacy of our method, we provide a new tool for investigating and controlling systems with unstable patterns, with potential implications for a wide range of scientific disciplines.
\end{abstract}

\maketitle

\begin{quotation}
Based on a Fourier basis, we characterize dynamically invariant subspaces and introduce convolution controls guided by symmetry groupoids. Our research highlights the effectiveness of symmetry groupoids in characterizing  spatio-temporal symmetries of patterns and their potential to facilitate pattern-selective feedback stabilization.
\end{quotation}

\section{\label{sec:introduction}Introduction}

Natural phenomena such as branching in plants \cite{MEI98}, leopard spots \cite{MUR88}, and sea shell \cite{MEI95}; as well as artificial phenomena like rotating spirals in the Belousov--Zhabotinsky reaction \cite{WIN84,ZHA91} all exhibit spatio-temporally ordered structures known as \emph{patterns}. Moreover, pattern formation is also observed in other fields, including fluid dynamics \cite{RAN82, CRA91}, solid-state physics \cite{HOY06}, and nonlinear optics \cite{LU96}. The diffusion-induced Turing instability reveals 
a universal mechanism for pattern formation \cite{TUR52}, making models in the form of reaction-diffusion partial differential equations (abbr. PDEs) highly suitable for describing and predicting a wide range of patterns \cite{FIE03, MUR03}.

From a mathematical perspective, a pattern can be defined as a representation of a solution to the underlying PDEs \cite{GOL03}. However, when the solution is unstable, it becomes vulnerable to errors in experimental or numerical simulations, resulting in a pattern that is rarely observable. Therefore, the control of unstable patterns has been studied extensively, both theoretically and in diverse applied contexts such as optics, lasers, or control of cardiac tissue \cite{HOE05, JEN98, KYR09, LEH11, OTT90, RAP99, SCH06, YAN06}. 

In this article, we explore a novel method to stabilize unstable solutions through the use of noninvasive feedback controls. More specifically, we introduce the concept of \emph{convolution controls}, which extends both the Pyragas control for ordinary differential equations (abbr. ODEs) \cite{PYR92} and the symmetry-breaking controls for PDEs \cite{SCH16, SCH22}. The specifics of convolution controls will be defined in the following sections. As a concrete application, we design convolution controls and exploit them to stabilize unstable equilibria (i.e., time-independent solutions) of the Chafee--Infante equation \cite{CHA74},
\begin{equation} \label{Chafee-Infante-equation}
\partial_t u = u_{xx} + \lambda u ( 1- u^2), 
\end{equation}
on the interval $x \in (0, \pi)$ equipped with Dirichlet boundary conditions 
\begin{equation} \label{boundary-condition}
u(t, 0) = u(t, \pi) = 0 \quad \mbox{for   } \, t \ge 0.
\end{equation} 
Here $\lambda > 0$ is a bifurcation parameter, which can be regarded as a scale on the length of the interval. Extensive research has been carried out on the Chafee--Infante equation\cite{HEN81,FIE96,FIE03}, and it is nowadays considered as a prototype for reaction-diffusion equations.

Our control strategy is motivated by the Pyragas control method, which is used to stabilize equilibria or periodic solutions of ODEs of the form $\dot{q}(t) = f(q(t))$, where $q(t) \in \mathbb{R}^n$. The Pyragas control method employs the following control system with time-delayed feedback:
\begin{equation} \label{pyragas-control}
\dot{q}(t) = f(q(t)) + b (q(t)  - q(t-\tau)).
\end{equation}
Here $\tau \ge 0$ is a prescribed time delay and $b \in \mathbb{R}^{n\times n}$ is called the feedback gain matrix. The \emph{Pyragas control term} 
\begin{equation} \label{Pyragas-control-term}
b (q(t) - q(t-\tau))
\end{equation}
is \emph{noninvasive} because it vanishes on equilibria and on any $\tau$-periodic solutions, while it possibly changes the local stability property of the targeted solution for a suitably chosen feedback gain $b$. The Pyragas control scheme is widely applicable because it it model independent and inexpensive to implement \cite{PYR06, SCH06, KOH09}. Various generalizations of the Pyragas control, in the setting of ODEs, have been established \cite{NAK98,SOC94, KIT95,SCH16a}.

\begin{figure}[t!]
\begin{center}
\includegraphics[width=0.49\textwidth]{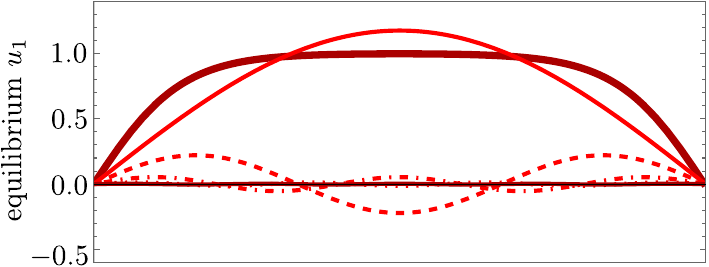}
\includegraphics[width=0.49\textwidth]{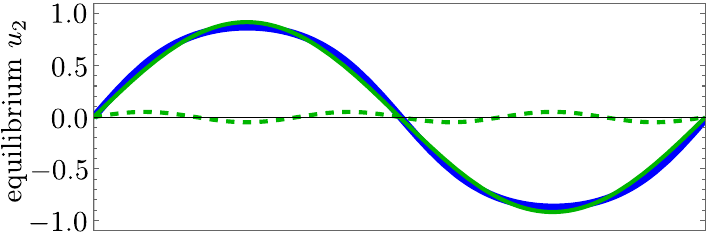}
\includegraphics[width=0.49\textwidth]{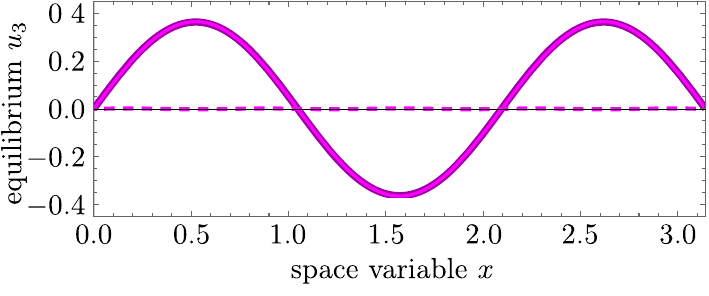}
\end{center}
\caption{
Equilibria of the Chafee--Infante equation \eqref{Chafee-Infante-equation}--\eqref{boundary-condition} for $\lambda=10$ together with their Fourier expansions for $k\leq 9$. The equilibrium $u_1$ (solid dark red curve) lies in the space $X_{1} = \left\{\sum_{\ell \in \mathbb{N}_{odd}} a_{1 \ell} \sin( \ell x)\right\}$, that is, it consists of all sine functions (red curves, the curve corresponding to $a_1 \sin(x)$ is solid, while the other curves are dashed, dot-dashed and dotted, respectively) that are even-symmetric with respect to $x=\frac{\pi}{2}$. In contrast, the equilibrium $u_2$ (solid dark blue curve) lying in the space $X_{2} = \left\{\sum_{\ell \in \mathbb{N}_{odd}} a_{2 \ell} \sin(2 \ell x)\right\}$ is odd-symmetric with respect to $x=\frac{\pi}{2}$ and locally even-symmetric with respect to $x=\frac{\pi}{4}, \frac{3\pi}{4}$. It consists of only those sine functions that obey those odd-even symmetries (sine curves in green, the curve corresponding to $a_2 \sin(2x)$ is solid, while the other curve corresponding to $a_6 \sin(6 x)$ is dashed). Last, the equilibrium $u_3$ (solid dark violet curve) is also locally even-symmetric with respect to $x=\frac{\pi}{6}, \frac{\pi}{2}, \frac{5\pi}{6}$, as well as locally odd-symmetric with respect to $x=\frac{\pi}{3}, \frac{2 \pi}{3}$. It therefore lies in a smaller subspace $X_{3} = \left\{\sum_{\ell \in \mathbb{N}_{odd}} a_{3 \ell} \sin(3 \ell x)\right\}$ than $X_1$ as compared to the equilibrium $u_1$, even though both equilibria would generally be described as ``even-symmetric'' with respect to $x = \frac{\pi}{2}$. Note that the curve corresponding to $a_3 \sin(3x)$ (solid magenta curve) lies on the equilibrium $u_3$ to a high precision, the curve corresponding to $a_9 \sin(9x)$ (dashed magenta curve) is almost identical zero.}
\label{general_example_sol-1}
\end{figure}

Motivated by the observation that symmetry-breaking bifurcations often trigger pattern formation in PDEs \cite{GOL03,HOY06}, Schneider extended the Pyragas control method to include symmetry-breaking controls \cite{SCH16}. These controls are designed by considering the spatio-temporal symmetries of the targeted solution rather than the full equivariance of the uncontrolled system. To illustrate this concept regarding controls, we consider the Chafee--Infante equation \eqref{Chafee-Infante-equation}--\eqref{boundary-condition}. The symmetry-breaking control system reads
\begin{equation} \label{control-Chafee-Infante-equation}
\partial_t u = u_{xx} + \lambda u (1-u^2) + b (u - (\pm 1) u(t-\tau, \pi - x)),
\end{equation}
where $u = u(t,x)$. The \emph{symmetry-breaking control term} 
\begin{equation} \label{symmetry-breaking-control-term}
b ( u - (\pm 1) u(t-\tau, \pi -x))
\end{equation}
includes a multiplicative factor $\pm 1$, a time delay $\tau \ge 0$, and a space reflection $x \mapsto \pi - x$. Hence it distinguishes between even-symmetric solutions (i.e., $u(t, \pi - x) = u(t,x)$) and odd-symmetric solutions (i.e., $u(t, \pi -x) = - u(t,x)$) with respect to $x = \frac{\pi}{2}$, but it does not distinguish other solutions. A theoretical limitation arises as a result, where symmetry-breaking controls are generally unable to stabilize all equilibria possessing the same reflection symmetry. In particular, it is impossible to select a specific equilibrium to be stabilized\cite{SCH16}; see Figure \ref{general_example_sol-1} for an illustration of two equilibria that possess the same reflection symmetry with respect to $x=\frac{\pi}{2}$. \textcolor{cyan}{}

In order to selectively stabilize \emph{all} unstable solutions, our article goes beyond the classical characterization of spatio-temporal symmetries of the targeted solution and proposes the design of more general control terms.

Our control design is motivated by the sifting property of the Dirac delta functional $\delta$, which can be expressed as 
\begin{equation} \label{shifting-property}
u(t, x)=(\delta\ast u)(t, x)=\int_0^\pi \delta(x - \xi) u(t, \xi) \, \mathrm{d}\xi,
\end{equation}
where $\ast$ denotes the standard convolution operation. This property implies that the control term $b(u - \delta \ast u)$ is noninvasive to all functions $u$. However, stabilization is bound to fail when the control term is identically zero. Therefore, it is crucial to identify suitable functionals such that the convolution control term is noninvasive exclusively on the targeted solution, while being invasive on other solutions.
This approach ensures that the targeted solution is selected and then stabilized.

It is important to emphasize that convolution controls represent a generalization of symmetry-breaking controls  
\eqref{symmetry-breaking-control-term}, because $u(t -\tau, \pi - x) = (h \ast u)(t- \tau, x)$ with $h(z) := \delta (\pi -z)$. 

Based on the Fourier basis associated with Dirichlet boundary conditions \eqref{boundary-condition} in Section \ref{sec:vertex-spaces} and a straightforward calculation in Section \ref{sec:vertex-groups}, we propose the following \emph{convolution control system}:
\begin{align} \label{convolution-Chafee-Infante-equation}
\partial_t u & = u_{xx} + \lambda u (1-u^2) + b (u - \mathcal{C}_{(h, \tau)}[u]),
\end{align}
where the \emph{convolution control operator} takes the form 
\begin{align} \label{convolution-control-operator}
\mathcal{C}_{(h, \tau)}[u](t,x)  = \int_0^\pi (h(x - \xi) - h(x+ \xi)) u(t - \tau, \xi) \, \mathrm{d}\xi. 
\end{align}
Here $h$ is a \emph{kernel functional} formally expressed by
\begin{equation} \label{formal-expression}
h(z) := \frac{1}{\pi} \sum_{m = 0}^\infty h_m \cos(m z), \quad z \in \mathbb{R}, \, h_m \in \mathbb{R}.
\end{equation}
The new control scheme will be \emph{pattern-selective} in the sense that it only preserves
the targeted equilibrium.

Our objective is therefore to identify the set $H_*$ of kernel functionals $h$ that allow stabilization of an unstable solution $u_*$. The key to stabilization lies  in the \emph{convolution control term} 
\begin{equation} \label{convolution-control-term}
b ( u - \mathcal{C}_{(h, \tau)}[u]),
\end{equation}
which should vanish solely on the targeted solution $u_*$, and not on any of the unstable and center eigenfunctions associated with $u_*$.
Therefore, it is essential to determine dynamically invariant subspaces $X_*$ that contain $u_*$, as well as to find suitable sets $H_*$. Achieving this requires more symmetries beyond equivariance \cite{CHO00,GOL88,GOL03} of the underlying model.

To demonstrate this requirement,
we consider the Chafee-Infante equation \eqref{Chafee-Infante-equation}--\eqref{boundary-condition}. In this case, the only equivariance presented is the space reflection, $x \mapsto \pi - x$, which represents an action of the group $\mathbb{Z}_2$. This equivariance only distinguishes solutions that are even-symmetric or odd-symmetric, and it is insufficient to stabilize all solutions with the same reflection symmetry. To overcome this limitation, we turn to the theory of \emph{symmetry groupoids}\cite{SCH22a} that allows us to incorporate additional symmetries that are necessary for stabilization.

The concept of groupoids extends the notion of groups \cite{BRO87,IBO19}. For the convenience of the reader, the definition of groupoids is given in Appendix \ref{defgroupoid}. Groupoids serve as a powerful instrument for understanding the symmetrical structure in different fields, such as quantum mechanics \cite{BAE01, CIA18}, string theory \cite{SAE13}, material science \cite{FRE13}, and general relativity \cite{BLO13}.

From an algebraic perspective, a groupoid is composed of a set of \emph{objects}, which, in our scenario, correspond to the dynamically invariant spaces $X_*$, and a set of \emph{morphisms}, which in our case, are represented by the kernel functionals that interconnect these objects. Unlike group elements, morphisms in a groupoid may not always have a defined composition, which offers more flexibility and intricacy in the symmetry structure. For the Chafee--Infante equation, groupoids enable us to identify symmetries on a functional level.

The concept of \emph{groupoid symmetries} (see Appendix \ref{yjxysymmetry}) has been developed by Schneider as a tool to characterize symmetries and equivariance in dynamical systems where the group theory may not be suitable \cite{SCH22a}, for instance, when symmetries are non-transitive or when various kinds of symmetries operate on distinct subspaces simultaneously.

This article is organized as follows: Section \ref{sec:bif} provides a brief explanation of the well-known bifurcation structure of Chafee--Infante equilibria. Section \ref{sec:vertex-spaces} delves into the characterization of dynamically invariant subspaces, referred to as vertex spaces, for any targeted equilibrium. In Section \ref{sec:vertex-groups}, sets of kernel functionals, called vertex isotropy groups, are characterized for stabilization through convolution controls. The proof of stabilization is presented in Section \ref{sec:stabilization}. Finally, Section \ref{sec:conclusion} offers a summary and explores several directions for future research.

To make the article more accessible, we have provided full mathematical details in the appendices. Specifically, Appendix \ref{sec:functionalsetting} introduces the functional setting and the notion of stability. Additionally, Appendix \ref{sec:groupoids} serves as a brief introduction to symmetry groupoids and equivaroid dynamical systems.

\section{Bifurcation structure of equilibria}\label{sec:bif}

To set the stage for our analysis and control in the upcoming sections, we give a concise explanation for the bifurcation structure of Chafee--Infante equilibria. To this end, we introduce the following notation: For each $s \ge 0$, $H_0^s$ denotes the Sobolev space of the functions $v :[0,\pi] \rightarrow \mathbb{R}$ such that $v(0) = v(\pi) = 0$ and the following norm is finite:
\begin{equation} \label{Sobolev-norm}
|v|_{H^s}:= \left(\sum_{k=0}^{\infty} (1 + k^{2})^s |\hat{v}_k|^2\right)^{\frac{1}{2}}.
\end{equation}
Here $\hat{v}_k := \frac{2}{\pi}\int_0^\pi v(x) \sin(kx) \, \mathrm{d}x$ denotes the $k$-th Fourier coefficient of $v$. The case $s = 0$ stands for the Lebesgue space $L^2$. 

\begin{figure}[t!]
\begin{center}
\includegraphics[width=0.45\textwidth]{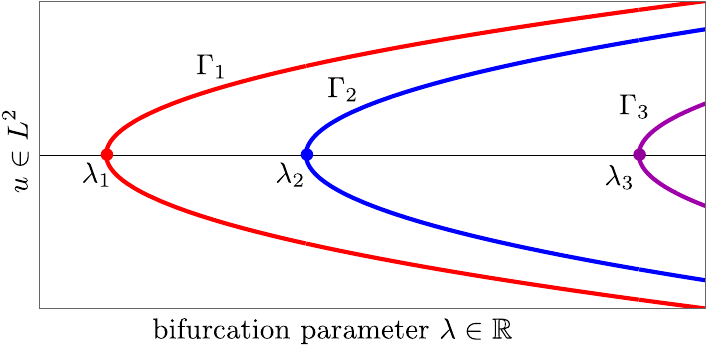}
\end{center}
\caption{The supercritical pitchfork bifurcation diagram of equilibria for the Chafee--Infante equation \eqref{Chafee-Infante-equation}--\eqref{boundary-condition}. The bifurcation curves are indexed by the unstable dimension $j-1 \in \mathbb{N}_0$ of equilibria. In particular, all equilibria on $\Gamma_j$ for $j \ge 2$ are unstable.}
\label{fig-bifurcation}
\end{figure}

It is known that the Chafee--Infante equation \eqref{Chafee-Infante-equation}--\eqref{boundary-condition} displays countably many bifurcation curves $\Gamma_j$ for $j \in \mathbb{N}$, which emanate from the \emph{trivial equilibrium} $u=0$. These bifurcation curves are of supercritical pitchfork type; see Ref. \cite{HEN81} Section 5.3 and Figure \ref{fig-bifurcation}. Each point $(\lambda, u)$ with $u \neq 0$ on a bifurcation curve corresponds to a nontrivial equilibrium. The curve $\Gamma_j$ bifurcates at 
\begin{equation} \label{bifurcation-value}
\lambda_j = j^2 \pi^2 \quad \mbox{for   } j \in \mathbb{N},
\end{equation}
corresponding to the $j$-th Dirichlet eigenvalue of the Laplace operator $-\partial_{xx}$ on the interval $[0,\pi]$. The $L^2$-normalized eigenfunction associated with $\lambda_j$ is $\frac{2}{\pi}\sin(jx)$.

According to the shooting argument used to solve the Sturm--Liouville eigenvalue problems (see Ref. \cite{DAI21a} Lemma 3.3 and Ref. \cite{DAI21} Theorem 3), every equilibrium $u_j(x)$ located on $\Gamma_j$ has $j-1$ simple zeros within the interval $(0,\pi)$, and its unstable dimension is also $j-1$. Therefore, all equilibria on $\Gamma_j$ with $j \ge 2$ are unstable and provide ideal candidates for stabilization.

Finally, using symmetry arguments from equivariant bifurcation theory (see Ref. \cite{DAI21a} Lemma 3.3), it is well established that for $j\in \mathbb{N}$ and $x\in [0,\pi]$, the relation 
\begin{equation}
    u_j(\pi - x) = (-1)^{j-1} u_j(x) 
\end{equation}
holds, indicating that every equilibrium has a reflection symmetry with respect to $x=\frac{\pi}{2}$.

In the following sections, we will unveil novel properties of the bifurcation curves. Concretely, we prove that these curves lie within distinctive dynamically invariant subspaces that we denote as \emph{vertex spaces}. Additionally, we will reveal specific actions of so-called \emph{vertex isotropy groups} that preserve each element of the vertex spaces. These vertex isotropy groups form part of the symmetry groupoid. All of these ingredients are crucial to our design of convolution controls.

\section{Dynamically invariant vertex spaces} \label{sec:vertex-spaces}

In order to characterize the dynamically invariant linear subspaces and consequently the symmetries of equilibria of the Chafee--Infante equation \eqref{Chafee-Infante-equation}--\eqref{boundary-condition}, we use the orthonormal Fourier basis of $L^2$,
\begin{equation} \label{orthonormal-basis}
\left\{\frac{2}{\pi}\sin(kx) : k \in \mathbb{N}\right\},
\end{equation}
whose elements are the $L^2$-normalized eigenfunctions of the Laplace operator $-\partial_{xx}$ with Dirichlet boundary conditions \eqref{boundary-condition}. 
Then, for each fixed $j \in \mathbb{N}$, we define the following $L^2$-subspaces, which we call \emph{vertex spaces} \cite{SCH22a}:
\begin{equation} \label{vertex-spaces}
X_{j} = \left\{\sum_{\ell \in \mathbb{N}_{odd}} a_{j \ell} \sin(j \ell x) \right\},
\end{equation}
where $\mathbb{N}_{odd}$ denotes the set of odd positive integers. See Figure \ref{general_example_sol-1} for a visualization of the equilibria lying in their respective vertex spaces.

\begin{remark}
We have the option to choose larger vertex spaces as an alternative, defined as
\begin{equation} \label{vertex-spaces2}
\tilde{X}_{j} = \left\{\sum_{\ell \in \mathbb{N}} \tilde{a}_{j \ell} \sin(j \ell x) \right\}.
\end{equation}
In this context, "larger" refers to the inclusion relationship $X_j \subset \tilde{X}_j$. While both options adequately describe the symmetry of the Chafee--Infante equation \eqref{Chafee-Infante-equation}--\eqref{boundary-condition}, in the context of control, we will observe that it is more advantageous to utilize the smaller vertex spaces $X_j$.
\end{remark}


The first important property of these vertex spaces $X_j$ is that they are \emph{dynamically invariant} for the Chafee--Infante equation \eqref{Chafee-Infante-equation}--\eqref{boundary-condition}, i.e., the nonlinearity 
\begin{equation} \label{nonlinearity}
\mathcal{F}(\lambda, u): = u_{xx} + \lambda u (1- u^2)
\end{equation}
maps a dense subspace $Y_j \subset X_j$ into $X_j$. To understand this in a rigorous way, we consider the functional setting $\mathcal{F} : H_0^2 \rightarrow L^2$; see Appendix \ref{sec:functionalsetting}. Accordingly, we define the $H^2$-subspace 
\begin{equation} \label{restriction}
Y_j = X_j\Big|_{H_0^2},
\end{equation}
i.e., the restriction of $X_j$ to $H_0^2$ equipped with the $H^2$-norm. Although the spaces $X_j$ and $Y_j$ differ by their norms, they consist of the same form for smooth functions and so they are identical from the symmetry perspective.

\begin{lemma}[Vertex spaces are dynamically invariant under $\mathcal{F}$] \label{lemma-invariance-F} 
The nonlinearity $\mathcal{F}(\lambda, u) =u_{xx}+ \lambda u(1-u^2)$ maps $Y_j$ into $X_j$. 
\end{lemma}
\begin{proof}[\textbf{Proof}]
Since $\{ \sin(kx): k \in \mathbb{N}\}$ is the set of Dirichlet eigenfunctions of $-\partial_{xx}$,  it suffices to show that the nonlinear part $u \mapsto u^3$ maps $Y_j$ into $X_j$. The Sobolev embedding theorem guarantees that every function in $H_0^2$ is continuous. So if $u \in Y_j \subset H_0^2$, then $u$ is continuous, and consequently, $u^3 \in L^2$. Finally, the trigonometric identity
\begin{equation} \label{trigonometric-identity-1}
4 \sin^3(\alpha) = 3 \sin(\alpha) - \sin(3 \alpha) \quad \mbox{for   } \alpha \in \mathbb{R}
\end{equation}
implies $u^3 \in X_j$.
\end{proof}

The next important property is that the bifurcation curve $\Gamma_j$ is contained in the vertex space $Y_j$.

\begin{lemma}[Nontrivial equilibria lie in vertex spaces] \label{lemma-space-for-equilibria} Let $u_j$ be a nontrivial equilibrium on $\Gamma_j$. Then $u_j \in Y_j$.
\end{lemma}
\begin{proof}[\textbf{Proof}]
Since $\mathcal{F}$ maps to $Y_j$ into $X_j$ by Lemma \ref{lemma-invariance-F}, we can restrict the Chafee--Infante equation \eqref{Chafee-Infante-equation}--\eqref{boundary-condition} to $Y_j$. At the bifurcation point $(\lambda_j, 0)$ of the bifurcation curve $\Gamma_j$, the linearized operator (i.e., the partial Fr\'{e}chet derivative) $\partial_u \mathcal{F}(\lambda_j, 0) : Y_j \rightarrow X_j$ defined by
\begin{equation} \label{partial-derivative}
\partial_u \mathcal{F}(\lambda_j, 0)[v] := v_{xx} + \lambda_j v
\end{equation}
has an eigenfunction $\sin(jx)$, which belongs to $Y_j$. Since $\lambda_j$ is a simple eigenvalue of $-\partial_{xx}$, by the theorem of bifurcation from simple eigenvalues \cite{CRA71}, nontrivial equilibria in $Y_j$ that bifurcate from $(\lambda_j, 0)$ exist and form a bifurcation curve $\tilde{\Gamma}_j \subset \mathbb{R} \times Y_j$. The uniqueness of bifurcation curves from simple eigenvalues implies $\Gamma_j = \tilde{\Gamma}_j$ and thereby $u_j \in Y_j$. 
\end{proof}

Notice that the same result as in Lemma \ref{lemma-space-for-equilibria} can be achieved by using the \emph{equivaroid branching lemma}, which is a broader variant of the classical equivariant branching lemma\cite{GOL03}; see Ref. \cite{SCH22a} for more information.

\section{Vertex isotropy groups and convolution controls} \label{sec:vertex-groups}

We begin by examining a targeted equilibrium $u_j$ on $\Gamma_j$ and we use the noninvasiveness condition to characterize sets of kernel functionals, which we refer to as vertex isotropy groups. After characterizing these groups, we derive an explicit expression for the convolution control operator. 

Recall that the targeted equilibrium $u_j$ belongs to the vertex space $Y_j$ by Lemma \ref{lemma-space-for-equilibria}. To stabilize $u_j$, our first step is to obtain kernel functionals $h$ such that the control term $b( u - \mathcal{C}_{(h, \tau)}[u])$ is noninvasive on $u_j$, i.e., 
\begin{equation} \label{noninvasive-to-uj}
u_j - \mathcal{C}_{(h, \tau)}[u_j] = 0.
\end{equation}
Since we aim to stabilize equilibria, the time delay $\tau \ge 0$ is irrelevant to \eqref{noninvasive-to-uj}. Hence for simplicity we take $\tau = 0$ and denote $\mathcal{C}_h := \mathcal{C}_{(h, 0)}$.

A sufficient condition for \eqref{noninvasive-to-uj} is to choose functionals $h$ such that the control term is noninvasive on all functions in $Y_j$. Equivalently, the control operator $\mathcal{C}_h$ acts as the identity operator on the targeted vertex space, i.e.,
\begin{equation} \label{Dirac-delta-property}
\mathcal{C}_{h}\Big|_{Y_j} = \mathcal{I}_j,
\end{equation}
where $\mathcal{I}_j : Y_j \rightarrow X_j$, $\mathcal{I}_j[v] = v$, is the identity operator. 

We collect all candidates to kernel functionals in the following \emph{vertex isotropy group}: For each fixed $j \in \mathbb{N}$, we define
\begin{equation} \label{vertex-group}
H_j := \left\{ \frac{1}{\pi} \sum_{m = 1}^\infty h_m \cos(m z) : h_{jm} =1 \mbox{   for   } m \in \mathbb{N}_{odd}\right\}
\end{equation}
and call each $h_m$ a \emph{control parameter}. Notice that $h(z) = \frac{1}{\pi}\sum_{m = 1}^\infty h_m \cos(m z)$ is not a function in any $L^p$-space, because $h_{jm} = 1$ for $m \in \mathbb{N}_{odd}$. Instead, the infinite series above is a \emph{formal expression} of the kernel functional $h$.

\begin{remark}
Had we chosen the larger vertex spaces $\tilde{X_j}$ defined in \eqref{vertex-spaces2} rather than $X_j$, it would have given us the smaller vertex isotropy groups
\begin{equation} \label{vertex-group2}
\tilde{H}_j := \left\{ \frac{1}{\pi} \sum_{m = 1}^\infty \tilde{h}_m \cos(m z) : \tilde{h}_{jm} =1 \mbox{   for   } m \in \mathbb{N}\right\},
\end{equation}
and so our choice of control parameters would have been limited. 
\end{remark}

It is important to recognize that the pair $(Y_j,H_j)$ constitutes a \emph{symmetry groupoid} in which the objects correspond to the vertex spaces $Y_j$ and the morphisms correspond to the vertex isotropy groups $H_j$. Notably, the symmetries within the vertex isotropy groups solely apply to the subspaces $Y_j$ and not to the full space $H_0^2$, rendering the standard group theory not applicable. Moreover, readers may observe that the vertex isotropy groups $H_j$ and $\tilde{H}_j$ do not contain the complete set of symmetries for the Chafee--Infante equation \eqref{Chafee-Infante-equation}--\eqref{boundary-condition}. The reason is that, for the sake of brevity, we have only included the part of the symmetry groupoid that is relevant to our control scheme. Furthermore, symmetry groupoids have applications in a broader range of contexts. The fundamental concepts and general definitions can be found in Appendix \ref{sec:groupoids}.

We are now ready to derive an explicit expression for the convolution control operator:
\begin{align} \label{concrete-convolution-control-operator-0}
\mathcal{C}_{h}[v](x) =  \int_0^\pi ( h(x-\xi) - h(x+\xi) ) v( \xi) \, \mathrm{d}\xi.
\end{align}
Note that this operator also  represents the action of the vertex isotropy group $H_j$ on the function space $H_0^2$: Given any kernel functional $h$  formally expressed by $h(z) = \frac{1}{\pi}\sum_{m = 1}^\infty h_m \cos(m z)$, the convolution control operator acts on functions $v(x) = \sum_{k = 1}^\infty a_k \sin(kx)$ in the full space $H_0^2$. Using the trigonometric identity
\begin{equation} \label{trigonometric-identity}
\cos(\alpha - \beta) - \cos(\alpha + \beta) = 2 \sin(\alpha)\sin(\beta) \quad \mbox{for   } \alpha, \beta \in \mathbb{R},
\end{equation}
and noting 
\begin{equation} \label{orthogonality}
\int_0^\pi \sin( m\xi) \sin (k \xi) \, \mathrm{d}\xi = \frac{\pi}{2}
 \delta_{m,k},
\end{equation}
where $\delta_{m,k}$ is the Kronecker delta defined by $\delta_{m, k} = 1$ if $m = k$ and $\delta_{m, k} = 0$ if $m \neq k$, from \eqref{concrete-convolution-control-operator-0}--\eqref{orthogonality} we obtain 
\begin{align} \label{concrete-convolution-control-operator-1}
\mathcal{C}_{h}[v](x) = \sum_{k = 1}^\infty h_k a_k \sin( kx).
\end{align}
Indeed, the convolution acts as a termwise multiplication on the Fourier coefficients of the function $v(x)$.

From the formula \eqref{concrete-convolution-control-operator-1}, it becomes apparent that the composition of two group elements from $H_j$ is simply a termwise multiplication of their respective control parameters. This composition yields another group element of $H_j$. If none of the control parameters  are zero, the group element is invertible, which is a relevant fact for describing its symmetry \cite{SCH22a}. Notice that invertibility is not necessary for our control scheme.

Therefore, we obtain an explicit expression of the convolution control term   
\begin{equation} \label{convolution-control-term-explicit}
b ( v - \mathcal{C}_h[v]) = b \sum_{k = 1}^\infty (1-h_k) a_k \sin(kx).
\end{equation}
The essence of the convolution control term is already evident from the formula \eqref{convolution-control-term-explicit}. The \emph{control parameters} $h_k$ of a kernel functional $h$ act as a filter for the corresponding Fourier mode $\sin(kx)$. Choosing $h_k = 1$ implies no control over the Fourier mode $\sin(kx)$, while choosing $h_k = -1$ (which we use for simplicity of analysis) means controlling the Fourier mode $\sin(kx)$. It is worth noting that the theoretical ability to adjust a countably infinite number of control parameters could be highly beneficial in practical applications of our control scheme.

The following lemma allows us to adopt the notion of stability and apply the principle of linearized stability in our functional setting; see Appendix \ref{sec:functionalsetting} for full details. 

\begin{lemma}[Boundedness and norm of control operators] \label{lemma-norm-control-operator}
Let $h$ be a kernel functional formally expressed by $h(z) = \frac{1}{\pi} \sum_{m = 0}^\infty h_m \cos(m z)$ with
\begin{equation} \label{smallness}
h_m \in \{-1, 1\} \quad \mbox{for   } m \in \mathbb{N}_0.
\end{equation} 
Then the convolution control operator $\mathcal{C}_{h}$ defined in \eqref{concrete-convolution-control-operator-0} is a linear bounded operator from $H_0^2$ to $H_0^2$. Moreover, its operator norm satisfies 
\begin{equation} \label{operator-norm}
\|\mathcal{C}_h\| := \sup_{v \in H_0^2 \setminus \{0\}} \frac{|\mathcal{C}_h[v] |_{H^2}}{|v|_{H^2}} \le 1.
\end{equation}
\end{lemma}
\begin{proof}[\textbf{Proof}]
The statements follow directly from the formula \eqref{concrete-convolution-control-operator-1} and the definition of Sobolev norm \eqref{Sobolev-norm}.
\end{proof}

Another result derived from the formula \eqref{concrete-convolution-control-operator-1} is the noninvasiveness of the convolution control operator on all functions $v \in Y_j$, provided that the kernel functionals $h$ are selected from the corresponding vertex isotropy group $H_j$.
 
\begin{lemma}[Noninvasiveness of the pair $(Y_j, H_j)$]\label{lemma-noninvasive-control} For each $h \in H_j$, the convolution control operator satisfies
\begin{equation} \label{dual-property}
\mathcal{C}_h[v] = v \quad \mbox{for   } v \in Y_j.
\end{equation}
In particular, the control term $b ( u - \mathcal{C}_h[u])$ is noninvasive to all nontrivial equilibria on the bifurcation curve $\Gamma_j$. 
\end{lemma}
\begin{proof}[\textbf{Proof}]
As we express $v(x) = \sum_{\ell \in \mathbb{N}_{odd}} a_{j\ell} \sin(j \ell x)$ and note $h_{j \ell} = 1$ for $\ell \in \mathbb{N}_{odd}$ due to $h \in H_j$, from \eqref{concrete-convolution-control-operator-1} we have 
\begin{align}
\begin{split}
\mathcal{C}_h[v](x) & =  \sum_{\ell \in \mathbb{N}_{odd}} h_{j \ell} a_{j \ell }  \sin(j \ell x)
\\&
=\sum_{\ell \in \mathbb{N}_{odd}} a_{j \ell} \sin(j \ell x)
\\&
= v(x),
\end{split}
\end{align}
which proves the claim.
\end{proof}

In order to achieve stabilization, it is crucial to ensure that the unstable and center eigenfunctions associated with $u_j$ do not belong to the vertex space $Y_j$, as the convolution control term \eqref{convolution-control-term-explicit} is noninvasive to $Y_j$. For the Chafee--Infante equation \eqref{Chafee-Infante-equation}--\eqref{boundary-condition}, it is well established that every equilibrium is hyperbolic, meaning that its associated center eigenspace is empty; see Ref. \cite{HEN81} Section 5.3. Therefore, it is only necessary to exclude the unstable eigenfunctions from $Y_j$.

\begin{lemma}[No instability in vertex spaces] \label{lemma-stabilization-must}
Let $u_j$ be a nontrivial equilibrium on $\Gamma_j$. Suppose that $v$ is an unstable eigenfunction associated with $u_j$. Then $v \notin Y_j$.
\end{lemma}
\begin{proof}[\textbf{Proof}]
Since the unstable dimension of $u_j$ is $j-1 \in \mathbb{N}_0$, by the shooting method (see the proof in Ref. \cite{DAI21a} Lemma 3.3), there exist $j-1$ unstable eigenfunctions $v_1, v_2, \cdots, v_{j-1}$ associated with $u_j$. All of these eigenfunctions $v_k$ with $k \in \{1,2,\cdots, j-1\}$ have $k-1$ simple zeros on $(0, \pi)$. 

Moreover, we also know that every function $w \in Y_j$ takes the form $w(x) = \sum_{\ell \in \mathbb{N}_{odd}} a_{j \ell} \sin(j \ell x)$ and $\sin(j x)$ has $j -1$ simple zeros on $(0, \pi)$. 
Using the zero number theorem \cite{ANG88}, we conclude that every nonzero function $w \in Y_j$ has at least $j-1$ simple zeros on $(0, \pi)$. This demonstrates that none of the unstable eigenfunctions of $u_j$ belong to $Y_j$.
\end{proof}

The remaining task is to prove that there exist some kernel functionals $h \in H_j$ such that the convolution control system 
\begin{equation} \label{concrete-convolution-control-system}
\partial_t u = u_{xx}  + \lambda u (1-u^2) + b (u - \mathcal{C}_{h}[u])
\end{equation}
stabilizes the targeted equilibria $u_j \in Y_j$.

\section{Stabilization} \label{sec:stabilization}

We complete the proof of stabilization for Chafee--Infante equilibria near bifurcation points. Through this section, let $(\lambda, u_j) \in \Gamma_j$ be a nontrivial Chafee--Infante equilibrium. Then the associated linearized operator $\mathcal{L}_{(b, h)}: H_0^2 \subset L^2 \rightarrow L^2$ with $(b, h) \in \mathbb{R} \times H_j$ reads
\begin{equation} \label{linear-control-system}
\mathcal{L}_{(b, h)}[v] := v_{xx} + \lambda (1 - 3 u_j^2) v + b (v - \mathcal{C}_h[v]).
\end{equation}
It is known that the spectrum of $\mathcal{L}_{(b, h)}$ consists solely of eigenvalues. Moreover, the so-called \emph{principle of linearized stability} holds, i.e., $u_j$ is locally exponentially stable (resp., unstable) if the real part of all eigenvalues of $\mathcal{L}_{(b, h)}$ is negative (resp., positive); see Appendix \ref{sec:functionalsetting} and Lemma \ref{lemma-principle-of-linearized-stability}. 

Our objective is to prove stabilization of all equilibria on $\Gamma_j$ that are near the bifurcation point $(\lambda_j, 0)$. To this end, by the upper semicontinuous dependence of spectra (see Ref. \cite{KAT95} Chapter 4, Remark 3.3), it suffices to consider $u_j = 0$ in \eqref{linear-control-system}, yielding the following simplified linearized operator: 
\begin{equation} \label{linear-control-system-near-zero}
\mathcal{L}_{(b, h)}[v] = v_{xx} + \lambda_j v + b (v - \mathcal{C}_h[v]).
\end{equation}
The fundamental characteristic of this operator is that the convolution controls preserve the set of eigenfunctions that we determined in  \eqref{orthonormal-basis}.

\begin{lemma}[Preservation of eigenfunctions under controls] \label{lemma-eigenfunctions-sharing} For all $(b, h) \in \mathbb{R} \times H_j$, the operator $\mathcal{L}_{(b, h)}$  possesses an identical set of eigenfunctions given by $\{\sin(kx) : k \in \mathbb{N}\}$.
\end{lemma}
\begin{proof}[\textbf{Proof}]
For each fixed $k \in \mathbb{N}$, we substitute $v_k(x) = \sin(k x)$ into the linearized operator \eqref{linear-control-system-near-zero}. Recalling that the convolution acts as a termwise multiplication (see \eqref{convolution-control-term-explicit}), we find 
\begin{align} \label{simple-eigen-relation}
\mathcal{L}_{(b,h)}[v_k](x) = ( -k^2 + \lambda_j + b (1 - h_k) ) v_k(x),
\end{align}
which implies that $v_k(x) = \sin(kx)$ is an eigenfunction of $\mathcal{L}_{(b,h)}$ for all $(b, h) \in \mathbb{R} \times H_j$. Given that $\{\sin(kx) : k \in \mathbb{N}\}$ forms a basis of $L^2$, its completeness as a basis implies that it serves as the set of eigenfunctions of $\mathcal{L}_{(b, h)}$ for all $(b, h) \in \mathbb{R} \times H_j$.
\end{proof}

We are ready to formulate our main result concerning stabilization near bifurcation points. Let $\lfloor a \rfloor$ denote the greatest integer smaller than $a \in \mathbb{R}$. 

\begin{thm}[Stabilization near bifurcation points] \label{theorem-near-0} Fix $j \in \mathbb{N}$. Consider the bifurcation point $(\lambda_j, 0) \in \Gamma_j$. Let 
\begin{equation} \label{choice-of-k}
\tilde{k}^2 := \lfloor \lambda_j \rfloor = \lfloor j^2\pi^2 \rfloor.
\end{equation} 
Suppose that the kernel functional $h(z) = \frac{1}{\pi} \sum_{m = 1}^\infty h_m \cos(mz)$ satisfies 
\begin{align} \label{finitary-near-0}
h_m = 
\left\{
\begin{array}{ll}
-1, \quad & \mbox{if   } m = 1, 2, \cdots, \tilde{k},
\\
1, & \mbox{if   } m \ge \tilde{k} + 1.
\end{array} \right.
\end{align}
Then the trivial equilibrium $u = 0$ becomes locally exponentially stable under the dynamics of the control system \eqref{concrete-convolution-control-system} with $\lambda = \lambda_j$ if we choose
\begin{align} \label{parameter-choice-near-0}
b < \frac{-\lambda_j}{2}.
\end{align}
Consequently, all nontrivial equilibria near the bifurcation point $(\lambda_j, 0)$ are locally exponentially stable under the dynamics of the control system \eqref{concrete-convolution-control-system}.
\end{thm}

\begin{proof}[\textbf{Proof}]
Let $\mu_k \in \mathbb{C}$ be an eigenvalue of $\mathcal{L}_{(b,h)}$ with $\sin(kx)$ as the associated eigenfunction. Then \eqref{simple-eigen-relation} implies 
\begin{equation} \label{characteristic-equation-near-0}
\mu_k = -k^2 + \lambda_j + b ( 1 - h_k).
\end{equation}
The conditions \eqref{choice-of-k}--\eqref{parameter-choice-near-0} imply $\mu_k < 0$ for all $k \in \mathbb{N}$. To see it, when $k = 1,2,\cdots, \tilde{k}$, we have
\begin{equation} \label{small-k-computation}
\mu_k < \lambda_j  + b( 1 - (-1)) < 0.
\end{equation}
When $k \ge \tilde{k} + 1$, we have $\mu_k = -k^2 + \lambda_j + b(1 - 1) < 0$. Hence $u = 0$ is locally exponentially stable according to the principle of linearized stability.
\end{proof}

\section{\label{sec:conclusion}Conclusion and Discussion}

In this article we delve into the possibility of utilizing noninvasive feedback controls to stabilize unstable solutions in reaction-diffusion PDEs. We introduce the concept of convolution controls, which not only extends the Pyragas control for ODEs and symmetry-breaking controls for PDEs, but also goes beyond the conventional characterization of spatio-temporal symmetries of solutions. Our design of the more general control terms is based on symmetry groupoids that explore symmetries of the underlying model from a functional point of view. Through the application of convolution controls, we successfully stabilize unstable equilibria of the Chafee--Infante equation \eqref{Chafee-Infante-equation}--\eqref{boundary-condition} near bifurcation points, demonstrating the effectiveness of our approach. Moreover, we are able to selectively stabilize equilibria of an arbitrarily high unstable dimension.

To apply the new control scheme in future applications, it is crucial to conduct a thorough analysis of the symmetry groupoid of the underlying model. Notably, the invariance property presented in Lemma \ref{lemma-invariance-F} remains applicable to analytic nonlinearities, wherein the nonlinear part is either an odd or even function. The proof utilizes the Taylor expansion and an analogous trigonometric identity for $\sin^n(\alpha)$, where $n \in \mathbb{N}$.

The effectiveness of convolution controls heavily relies on minimizing the size of vertex spaces $Y_j$ while simultaneously maximizing the size of vertex isotropy groups $H_j$. For example, merely selecting the vertex space $Y_1$ of odd functions does not guarantee stabilization of any equilibrium also lying in a smaller vertex space, as demonstrated in Ref. \cite{SCH16}. It is important to note that the case study of the Chafee--Infante equation presented in this article has the potential to be extended to any other dynamical system as long as careful consideration is given to choose vertex spaces to ensure stabilization by convolution controls.

Regarding the practical implementation of our convolution controls, notice that distributed controls, which involve convolutions or a Fourier basis, have been utilized in engineering for a considerable period \cite{BAM02, MAR96}. Therefore, the novelty is not the convolution process itself, but rather our method of employing convolution as a means of pattern selection by exploiting symmetries at a functional level. 

\begin{acknowledgments}
I.S. has been supported by the Deutsche Forschungsgemeinschaft, SFB 910, Project A4 “Spatio-Temporal Patterns: Control, Delays, and Design”. J.-Y. D. has been supported by MOST grant number 110-2115-M-005-008-MY3. We would like to express our gratitude to all the members of SFB 910 for their valuable contributions and fruitful discussions during the twelve years of funding. We extend our special thanks to Sabine Klapp and Eckehard Schöll. Additionally, we would like to acknowledge Bernold Fiedler, Alejandro L\'{o}pez Nieto, and Babette de Wolff for their insightful and inspiring discussions.
\end{acknowledgments}

\section*{Data Availability Statement}

Data sharing is not applicable to this article as no new data were created or analyzed in this study.

\appendix

\section{\label{sec:functionalsetting}Functional setting and notion of stability}

The Laplace operator $-\partial_{xx}: H_0^2 \subset L^2 \rightarrow L^2$ with $H_0^2$ as the domain generates a linear analytic semiflow on $L^2$. Moreover, the elliptic operator $\mathcal{F}: \mathbb{R} \times H_0^2 \rightarrow L^2$, defined by
\begin{equation}
\mathcal{F}(\lambda, u) = u_{xx} + \lambda u(1-u^2)
\end{equation}
generates a nonlinear $C^0$-semiflow $\{\mathcal{S}(t)\}_{t\geq 0}$ on any fractional space $H^{2\gamma}$ with $\gamma > 1/2$; see Ref. \cite{HEN81} Theorem 3.3.3.

Let $u_* \in H_0^2$ be an equilibrium of the Chafee--Infante equation \eqref{Chafee-Infante-equation}--\eqref{boundary-condition} at a parameter value $\lambda = \lambda_*$, i.e., $\mathcal{F}(\lambda_*, u_*) = 0$. 

\begin{defi}[Notion of stability]
The equilibrium $u_*$ is called \emph{locally exponentially stable} if for any $\varepsilon > 0$, there exist  $r > 0$ and $\alpha > 0$ such that for every $u_0 \in H^{2\gamma}$ satisfying $|u_* - u_0|_{H^{2 \gamma}} < r$, we have $|\mathcal{S}(t)[u_0] - u_*|_{H^{2\gamma}} < \varepsilon$ for $t \ge 0$ and $\lim_{t \rightarrow \infty} e^{-\alpha t}|\mathcal{S}(t)[u_0] - u_*|_{H^{2\gamma}} = 0$.
\end{defi}

By using the variation-of-constants formula, we have the following criterion for the local exponential stability; see Ref. \cite{HEN81} Chapter 5.
\begin{lemma}[Principle of linearized stability] \label{lemma-principle-of-linearized-stability}
The equilibrium $u_*$ is locally exponentially stable (resp., unstable) under the dynamics of the semiflow $\{\mathcal{S}(t)\}_{t \ge 0}$ if the spectrum of the linearized operator $\mathcal{L}_*: H_0^2 \subset L^2 \rightarrow L^2$, defined by
\begin{equation} \label{linearized-operator}
\mathcal{L}_*{[v]} := v_{xx} + \lambda_* (1 - 3 u_*^2) v,
\end{equation}
is contained in the left half-plane $\left\{z \in \mathbb{C} : \mathrm{Re}(z) < 0\right\}$ (resp., right half-plane $\left\{z \in \mathbb{C} : \mathrm{Re}(z) > 0\right\}$).
\end{lemma}

Suppose that the control system \eqref{concrete-convolution-control-system} is noninvasive on the equilibrium $u_*$, i.e., $u_*$ is also an equilibrium of \eqref{concrete-convolution-control-system}. Then the associated linearized operator $\mathcal{L}_{(b, h)} : H_0^2 \subset L^2 \rightarrow L^2$ reads
\begin{equation} \label{linear-control-system-appendix}
\mathcal{L}_{(b, h)}[v] := \mathcal{L}_*[v] + b ( v - \mathcal{C}_h[v])
\end{equation}
Suppose that the assumption \eqref{smallness} in Lemma \ref{lemma-norm-control-operator} holds. Then $\mathcal{C}_h : H_0^2 \rightarrow H_0^2$ is bounded. By the elliptic regularity theory and Sobolev embedding theorem \cite{HEN81}, $\mathcal{L}_{(b, h)}$ has compact resolvent, and thus its $L^2$-spectrum consists solely of eigenvalues. Moreover, the principle of linearized stability in Lemma \ref{lemma-principle-of-linearized-stability} also holds in the control setting \eqref{linear-control-system-appendix}.

\section{\label{sec:groupoids} Symmetry groupoids}
In what follows, we first give the algebraic definition of groupoids and their connection to symmetry. This appendix closely follows Ref. \cite{SCH22a}.

\begin{defi}[Groupoids \cite{SCH22a, IBO19, OLV15, WEI96}]\label{defgroupoid}
Let $B$ be a set. A \emph{groupoid} is a set $\Gamma$ of morphisms $\gamma: B \rightarrow B$, $\gamma \in \Gamma$, supplemented with the following maps:
\begin{itemize}
\item a surjective \emph{source map} $s\colon \Gamma \rightarrow B$, $\gamma \mapsto s(\gamma)$,
\item a surjective \emph{target map} $t\colon \Gamma \rightarrow B$, $\gamma \mapsto t(\gamma)$,
\item a \emph{partial binary operation} $\circ$ defined on the set of composable morphisms $\Gamma \star \Gamma:=\left\{ (\gamma_2,\gamma_1) \in \Gamma \times \Gamma : t (\gamma_1)=s(\gamma_2) \right\}$:
\begin{equation}
\begin{split}
\circ\colon \quad  \Gamma \star \Gamma & \rightarrow \Gamma\\
(\gamma_2,\gamma_1)& \mapsto \gamma_2\circ \gamma_1,
\end{split}
\end{equation} 
\item an injective \emph{identity map} $e \colon B\rightarrow \Gamma$, $b\mapsto e(b)=: e_b$,
\end{itemize}
which satisfy the following properties:
\begin{enumerate}
\item the partial binary operation is associative, that is, for all $(\gamma_3,\gamma_2)$, $(\gamma_2,\gamma_1)\in \Gamma\star \Gamma$, the identity  $(\gamma_3 \circ \gamma_2)\circ \gamma_1=\gamma_3 \circ (\gamma_2  \circ \gamma_1)$ holds;
\item the identity map defines a family of identity morphisms in the following sense:
\begin{enumerate}
\item for all $b \in B$: $s(e_b)=t(e_b)=b$,
\item for all $\gamma$ such that $s(\gamma)=b$:  $\gamma \circ e_{b}=\gamma$,
\item for all $\gamma$ such that $t(\gamma)=b$:  $e_b \circ \gamma=\gamma$;
\end{enumerate}
\item each morphism $\gamma \in \Gamma$ has a two-sided inverse $\gamma^{-1}\in \Gamma$ such that \begin{equation}
s(\gamma) = t(\gamma^{-1}), \quad t(\gamma) = s(\gamma^{-1}), \quad \mathrm{and}
\end{equation} 
\begin{equation}
\gamma^{-1} \circ \gamma= e_{s(\gamma)}, \quad \gamma \circ\gamma^{-1}= e_{t(\gamma)}.
\end{equation}
\end{enumerate} 
We denote such a groupoid by $(\Gamma\rightrightarrows B)$. The set $B$ is called the \emph{base}, and its elements are called \emph{objects}. Moreover, we call $s(\gamma)$ the \emph{source} of the morphism $\gamma$, and $t(\gamma)$ its \emph{target}.
\end{defi}

The abstract algebraic concept of groupoids and a concrete dynamical system are connected via the following definition.

\begin{defi}[Groupoid symmetry for $C^0$-semiflows \cite{SCH22a}]\label{yjxysymmetry}
Consider a $C^0$-semiflow $\{\mathcal{S}(t)\}_{t \geq 0 }$ on a Banach space $X$. Let $\{X_j\}_{j \in B}$ be an indexed family of linear closed subspaces (``vertex spaces'') of $X$ such that $\mathcal{S}(t)X_j \subseteq X_j$ for all $t \geq 0$ and $j \in B$. Let $(\Gamma \rightrightarrows B)$ be a groupoid $\Gamma$ over the base $B$. We say that $\gamma \colon j\rightarrow k$, $\gamma \in \Gamma$, is a \emph{$(X_j,X_k)$-symmetry} of the semiflow $\{\mathcal{S}(t)\}_{t \geq 0 }$ if the following holds: There exists a representation $\rho_X$ of the groupoid $(\Gamma\rightrightarrows B)$ on the space $X$ such that
\begin{equation}\label{liniso_inf}
\rho_X(\gamma)\colon  X \rightarrow X \  \mathrm{with} \ \rho_X(\gamma) X_j=X_k 
\end{equation}
and 
\begin{equation}
\mathcal{S}(t)\rho_X(\gamma) y_0=\rho_X(\gamma) \mathcal{S}(t) y_0
\end{equation}
holds for all $y_0 \in X_j$ and $t \geq 0$.
\end{defi}	
The largest groupoid whose non-identity morphisms all represent nontrivial $(X_j,X_k)$-symmetries of a semiflow is called the \emph{symmetry groupoid}.

Consider a vertex space $X_j$, where $G_{jk}$ represents the set of $(X_j,X_k)$-symmetries. It is worth noting that the set $G_j:=G_{jj}$ of $(X_j,X_j)$-symmetries forms a group, referred to as the \emph{vertex symmetry group} of $X_j$.  Importantly, the vertex symmetry group $G_j$ preserves the vertex space $X_j$ as a set, meaning that $\gamma X_j=X_j$ for all $\gamma \in G_j$.
In addition to the vertex symmetry group $G_j$, there exists a subset within it that preserves the vertex space $X_j$ pointwise. This subset is referred to as the \emph{vertex isotropy group}.

Regarding the Chafee--Infante equation \eqref{Chafee-Infante-equation}--\eqref{boundary-condition}, it is worth noting that the vertex spaces are given by the spaces $X_j$ defined in \eqref{vertex-spaces}, $\tilde{X_j}$ defined in \eqref{vertex-spaces2}, and $X_0:=\{0\}$. The respective vertex symmetry groups are given by
\begin{equation} \label{vertex-group3}
G_j := \left\{ \pm \frac{1}{\pi} \sum_{m = 1}^\infty h_m \cos(m z) : h_{jm} =1 \mbox{   for   } m \in \mathbb{N}_{odd}\right\},
\end{equation}
\begin{equation} \label{vertex-group4}
\tilde{G}_j := \left\{ \pm \frac{1}{\pi} \sum_{m = 1}^\infty h_m \cos(m z) : h_{jm} =1 \mbox{   for   } m \in \mathbb{N}\right\},
\end{equation}
\begin{equation} \label{vertex-group5}
G_0 := \left\{ \pm \frac{1}{\pi} \sum_{m = 1}^\infty h_m \cos(m z)\right\}.
\end{equation}
The vertex isotropy groups $H_j$, $\tilde{H}_{j}$, and $H_0$ are defined by taking the positive sign in respective vertex symmetry groups.


\bibliography{aipsamp.bib}

\end{document}